\documentclass{article}
\usepackage{spconf,amsmath, amsthm}

\usepackage[T1]{fontenc}
\usepackage{color}
\usepackage{psfrag}
\usepackage[dvips]{graphicx}
\usepackage{tikz}
\usepackage{pgfplots}
\pgfplotsset{compat=1.12}
\usepackage{enumerate}
\usepackage{rotating}
\usepackage{srcltx}
\usepackage{booktabs}
\usepackage{multicol}
\usepackage{enumitem}
\usepackage{amsthm}

\usepackage{mathtools}
\usetikzlibrary{arrows,decorations.pathmorphing,backgrounds,positioning,fit,matrix}

\usepackage[noadjust]{cite}
\usepackage{amssymb}
\usepackage{dsfont}    
\usepackage{mathtools} 
\usepackage{stackengine}
\usepackage[toc,page]{appendix}
\usepackage{bbm}

\usepackage{cuted}
\usepackage{multicol,lipsum}
\usepackage{cuted}

\newtheorem*{corollary*}{Corollary}
\newtheorem{prop}{Proposition}
\newtheorem{lem}{Lemma}
\newtheorem{theorem}{Theorem}
\newtheorem*{theorem*}{Theorem}
\newtheorem{remark}{Remark}

\usepackage{algorithmic}
\usepackage{array}

\begin{document}
	
	\maketitle
	
	\begin{abstract}
 This work establishes regularity conditions for consistency and asymptotic normality of the multiple parameter maximum likelihood estimator (MLE)  from censored data, where the censoring mechanism is in the form of $1$-bit measurements.
 The underlying distribution of the uncensored data is assumed to belong to the exponential family, with natural parameters expressed as a linear combination of the predictors, known as generalized linear model (GLM). As part of the analysis, the Fisher information matrix is also derived for both censored and uncensored data, which helps to quantify the impact of censoring and assess the performance of the MLE. The choice of GLM allows one to consider a variety of practical examples where 1-bit estimation is of interest. In particular, it is shown how the derived results can be used to analyze two practically relevant scenarios: the Gaussian model with both unknown mean and variance, and the Poisson model with an unknown mean.

	\end{abstract}
	
	\keywords{$1$-bit measurements, censored data, exponential family, maximum likelihood, generalized linear model, Fisher information, radar, laser communications, Gaussian noise, Poisson noise.}
	
	\section{Introduction} 
	\label{sec:Intro}
In many engineering and scientific fields, parameter estimation from censored data is a critical challenge due to the inherently limited information available from measurements or observations. Such scenarios frequently arise in many applications, such as: 
(i) {\em sensor fusion}, where power and bandwidth constraints force sensors to quantize observations to a single bit, effectively censoring the data~\cite{kumar2023communication,fang2013}; 
(ii) {\em radar and wireless systems}, where the use of 1-bit analog-to-digital converters at the receiver end offers a power-efficient solution for future wireless systems that handle large signal bandwidths or numerous radio frequency chains~\cite{studer2016quantized}; 
and (iii) {\em survival analysis}, where the event of interest (e.g., system failure or patient mortality) may not be observed within the study period~\cite{jewell2003current}. 
	
In several applications where censoring mechanisms are applied or advantageous, the underlying distributions often belong to the \emph{exponential family}. For instance, Gaussian noise is widely used in numerous communication systems. Similarly, Poisson noise, commonly called shot noise, is frequently encountered in photon-limited imaging systems, optical communication, and various counting processes~\cite{yang2011bits}. The exponential distribution is broadly used in survival analysis due to its suitability for modeling time-to-event data~\cite{alzeley2021}.

While the maximum likelihood estimator (MLE) has been widely employed for parameter estimation from censored data, previous studies have often been limited in scope, focusing on \emph{specific distributions} tailored to particular applications and estimating a \emph{single unknown parameter}. Furthermore, there has been a lack of theoretical analysis regarding the (necessary or sufficient) \emph{regularity conditions} to ensure \emph{consistency} and \emph{asymptotic normality}.
This work aims to address this gap by providing a comprehensive study of the asymptotic behavior of the MLE under a rather general model of censored data from an exponential family, where the censoring mechanism is in the form of $1$-bit measurements.

In the signal processing domain censored, or 1-bit data, is often studied under an additive noise model (i.e., ordinary linear regression model).  In this work, we take a more general and unifying approach by adopting a generalized linear model (GLM)~\cite{myers2012}.  The key idea behind the adopted GLM is that it assumes that the natural parameters of the exponential family distribution are described by a linear combination of the predictors, hence making it a versatile tool in these contexts. 

    
    \subsection{Relevant literature}
Parameter estimation with $1$-bit data has been widely studied, and a comprehensive review is beyond this paper's scope. Instead, we focus on highlighting key relevant works.

Estimating the mean from 1-bit quantized signals has been extensively explored in various domains, including sensor fusion, radar, and communication systems. Previous studies, such as those in~\cite{fang2013,cheng2019target}, employed the MLE to estimate the unknown mean under the assumption of Gaussian noise. Notably, these works primarily focused on single parameter estimation where the variance is known. 
Recently, the author of~\cite{kumar2025onebitdistributedmeanestimation} proposed and analyzed a new estimator for unknown mean and variance where samples are, i.i.d from a scale-location family of distributions. 
The GLM approach taken in this paper allows one to also study models where multiple parameters are unknown.
A rather general and complete answer to the estimation of the mean from censored data under an additive noise model was provided in~\cite{kipnis2022}. In particular, the authors of~\cite{kipnis2022} derived asymptotically optimal estimators and the rate of convergence for a censored model with additive noise coming from a log-concave distribution. Additionally, the authors of~\cite{kipnis2022} provided these results under various cooperation models between the 1-bit compressors.
 
The focus on more general parameter estimation that extends beyond mean estimation has received less attention.
The authors of~\cite{dabeer2008} provided a more in-depth treatment of the MLE for an additive noise model with multiple parameters, allowing correlation across noise  samples. 
In~\cite{stoica2021cramer}, the authors
analyzed the Fisher
information matrix (FIM) for censored data in the context of linear models across arbitrary distributions, providing important theoretical foundations for parameter estimation in these settings. Estimation of parameters in auto-regressive processes with additive noise from censored data was considered in~\cite{1bitAutoRegressive, krishnamurthy1996asymptotic}. 
Kalman and particle filter methods for estimating signals from censored data were explored in works like~\cite{jeon2019reinforcement, duan2008state}. Recent approaches, such as~\cite{khobahi2019, jeon2019reinforcement}, employed deep neural networks and reinforcement learning under Gaussian noise assumptions.

\subsection{Our Contributions and Paper Outline}
The paper contribution and outline are as follows:
\begin{enumerate}[leftmargin=*]
\setlength\itemsep{0.01em}
 \item  Section~\ref{sec:ProbForm} presents the censored GLM and  the 
 MLE. 
\item Section~\ref{sec:Main Result} presents our main results. 
Proposition~\ref{prop:FIM}, in Section~\ref{subsec:Fisher Information}, characterizes the structure of the 
FIM
under the censored GLM. Theorem~\ref{thm:MainThm}, in Section~\ref{sec:cons_normality}, provides a set of conditions, which can be considered mild, that guarantee the consistency and asymptotic normality of the MLE.
\item Section~\ref{sec:Examples} demonstrates the versatility of our approach by focusing on several practically relevant examples. We analyze the optimal asymptotic performance of the MLE for two common noise distributions: Gaussian and Poisson. Notably, we analyze the case of unknown mean and variance for Gaussian noise, which is practically important (e.g., 1-bit radar) but has received limited attention. For both the Gaussian and Poisson models, we establish simple conditions for consistency and asymptotic normality.
\end{enumerate}

	\noindent  {\bf Notation.} For any $k \in \mathbb{N}$, we define $[1:k] = \{1,2,\ldots, k\}$; 
 logarithms are in base $\rm{e}$; $\mathbf{0}_k$ is the column vector of dimension $k$ of all zeros; $\|\mathbf{x}\|$ is the L2 norm of $\mathbf{x}$ and $\|\mathbf{x}\|_\infty$ is the L-infinity norm of $\mathbf{x}$.
 For two nonnegative definite matrices $\mathbf{A}$ and $\mathbf{B}$ of equal size, we let $\mathbf{A} \succeq \mathbf{B}$ indicate that $\mathbf{A} - \mathbf{B}$ is nonnegative definite.
 We use the notion of weak consistency.
	%
	%
	\section{Problem Formulation}
	\label{sec:ProbForm}
	Let $\mathbf{X} = \begin{bmatrix} X_1 & X_2 & \ldots & X_n \end{bmatrix}^T$ be an $n$-dimensional random vector with independent components $X_i \in \mathcal{X}$, such that 
 its density is $p_{\mathbf{X}}(\mathbf{x} ;\boldsymbol{\theta} ) = \prod_{i=1}^n p_{X_i}(x_i; \boldsymbol{\theta}),$ where
	\begin{equation}
		\label{eq:PDFExpFamily}
		p_{X_i}(x;\boldsymbol{\theta} ) = h(x) \exp (\langle \boldsymbol{\eta}_{i,\boldsymbol{\theta}}, \mathbf{T}_x \rangle - \phi(\boldsymbol{\eta}_{i,\boldsymbol{\theta}})), 
	\end{equation}	
	that is, each $X_i$ has a density that belongs to the \emph{exponential family}.
	With reference to~\eqref{eq:PDFExpFamily}, we have that:
	(i) $h(x): \mathcal{X} \to (0,+\infty)$ is the base measure;
	(ii) $\mathbf{T}_x: \mathcal{X} \to \mathbb{R}^d$ is the $d$-dimensional sufficient statistic vector;
	(iii) $\phi: \mathcal{X} \to \mathbb{R}$ is the log-partition function;
	and (iv) $\boldsymbol{\eta}_{i,\boldsymbol{\theta}}$ is the 
	$d$-dimensional  natural parameter vector, where in our model $\boldsymbol{\eta}_{i,\boldsymbol{\theta}} = \mathbf{V}_i\boldsymbol{\theta}$ with $\boldsymbol{\theta} \in \Theta$ being a $k$-dimensional vector that is common to $\boldsymbol{\eta}_{i,\boldsymbol{\theta}}$, for all $i \in [1:n]$ and needs to be estimated, and $\mathbf{V}_i$ is a known $d \times k$  matrix. 
 We will refer to this setting as a GLM.

	We assume that  $\mathbf{X}$ is not observed directly, and instead the measurements are given by  a $1$-bit censoring mechanism:
	\begin{equation} \label{eq:BinaryRV}
		B_i = \begin{cases}1 & \text { if } X_i \leq \tau_i,
			\\ -1 & \text { if } X_i > \tau_i,
		\end{cases}
	\end{equation}
	where $\tau_i$ for $i \in [1:n]$ is a fixed threshold and an interior point of $\mathcal{X}$.
	Therefore, for $b \in \{-1,1\}$, we have that\footnote{This should be understood in the sense of Lebesgue integral where $\rm{d}x$ is some appropriate dominating measure. In particular, for discrete distributions, it should be understood as  a summation.}
	\begin{equation}
		\label{eq:PMFBi}
		P_{B_i}(b; \boldsymbol{\theta}) = \Pr(B_i=b)
		= \int \limits_{\mathcal{X}(b)} p_{X_i}(x;\boldsymbol{\theta} ) \ {\rm{d}}x,
	\end{equation}
	where $\mathcal{X}(b)$ is a subset of $\mathcal{X}$ determined by $b$; specifically, $\mathcal{X}(1) = \{x\in\mathcal{X}:x\leq \tau_i\}$ and $\mathcal{X}(-1) = \{x\in\mathcal{X}:x > \tau_i\}$.
	
	In many applications,  as discussed in Section~\ref{sec:Intro}, one can only observe the censored data $\left\{ b_i \right \}_{i=1}^n$, which are realizations of $\left\{ B_i \right \}_{i=1}^n$. From these observed data samples we seek to estimate the $k$-dimensional vector $\boldsymbol{\theta}$ using the MLE method. In particular, the MLE is defined as follows\footnote{We note that the MLE may not be unique. In this case, we randomly select one of these possible choices.}, 
	\begin{equation}
		\label{eq:MLE}
		\boldsymbol{\hat{\theta}}_n 
		= \underset{\boldsymbol{\theta} \in \Theta}{\arg \max } \ \ell_n \left (\boldsymbol{\theta} ; \left \{ b_i \right \}_{i=1}^n \right ) ,
	\end{equation}		
	where $\ell_n \left (\boldsymbol{\theta} ; \left \{ b_i \right \}_{i=1}^n \right )  = \sum_{i=1}^n \log \left( P_{B_i}(b_i;\boldsymbol{\theta}) \right )$ is the log-likelihood function. 
	By using~\eqref{eq:PMFBi}, it is not difficult to show that (see Appendix~\ref{app:ProofSimplifLogLik} for the detailed computation),
	\begin{equation}
		\ell_n \left (\boldsymbol{\theta};\{b_i\}_{i=1}^n \right ) = \sum \limits_{i=1}^n 	\left[\phi(\boldsymbol{\eta}_{i,\boldsymbol{\theta}};b_i) - \phi(\boldsymbol{\eta}_{i,\boldsymbol{\theta}})\right],
	\label{eq:LogLikelihoodFunction}
	\end{equation}
	where $\phi(\boldsymbol{\eta}_{i,\boldsymbol{\theta}};b_i)$ is the log-partition function of  $p_{X_i | B_i}(x|b_i)$.  

	\section{Main Result} \label{sec:Main Result}
 Here, we present our main result.
	We start by characterizing the FIM, which plays a significant role in our analysis.
	
	\subsection{Fisher Information} \label{subsec:Fisher Information}
	The FIM is important for many reasons.
	First, it characterizes the asymptotic variance of the MLE, as we show in Section~\ref{sec:cons_normality}; second, via Cram\'er-Rao~\cite{cramer1999mathematical}, it can be used to provide a fundamental lower bound on some performance metrics, such as the mean square error and variance. 
    Stoica et al.~\cite{stoica2021cramer} analyzed the FIM for censored data in the context of linear models across arbitrary distributions. 
    In the following proposition, we extend the characterization of the FIM for GLMs, as discussed in Section~\ref{sec:ProbForm}.
	
	\begin{prop}
		\label{prop:FIM}
		The FIM of estimating the true parameter vector $\boldsymbol{\theta}_0$ from the censored data $\{B_i\}_{i=1}^n$ is given by 
		\begin{equation} \label{eq:FIM_censored}
			\mathbf{J}_n  = \sum_{i=1}^n \mathbf{V}_i^T \mathsf{Cov} \left( \mathbb{E} \left [ \mathbf{T}_{X_i}|B_i\right ] \right ) \mathbf{V}_i.
		\end{equation}
	\end{prop}
	\begin{proof}
		The proof is provided in Appendix~\ref{app:FIM}.
	\end{proof}
It is interesting to compare 
the FIM in~\eqref{eq:FIM_censored} to the FIM  of estimating $\boldsymbol{\theta}_0$ from the uncensored data $\{X_i\}_{i=1}^n$, which following similar computations as in Appendix~\ref{app:FIM} is given by
	\begin{equation} \label{eq:FIM_uncensored}
		\mathbf{I}_n = \sum_{i=1}^n \mathbf{V}_i^T \mathsf{Cov} \left( \mathbf{T}_{X_i} \right ) \mathbf{V}_i.
	\end{equation}
 One can establish an inequality between the two  FIMs. Specifically, $\mathbf{I}_n \succeq \mathbf{J}_n$, which follows from the data-processing  inequality of the Fisher information \cite{ZamirIT1998} using the Markov chain  $\boldsymbol{\theta}_0 \to X_i \to B_i$ for $i \in [1:n]$. This suggests a loss of information, as expected, when using censored data. The magnitude of this loss will depend on the considered distribution and choice of threshold $\tau_i$ for $i \in [1:n]$. Note that choosing a set of optimal $\tau_i$'s will help to maximize the FIM; however, this choice often relies on unknown parameters that need to be estimated. 
 In Section~\ref{sec:Examples}, we will provide
 a choice of the optimal thresholds for the Gaussian distribution.
 

\subsection{Consistency and Asymptotic Normality of the MLE} 
	\label{sec:cons_normality}
	The next theorem, whose proof is provided in Appendix~\ref{app:MainTheorem}, states the regularity conditions under which the MLE in~\eqref{eq:MLE}
 is \emph{consistent} and \emph{asymptotically normal}.
	\begin{theorem}
		\label{thm:MainThm}
		Provided that the following conditions hold:
		\begin{enumerate}
\setlength\itemsep{0.01em}
			\item \label{con:1 Thm}$\lim \limits_{n \to \infty} \ \max \limits_{i \in [1:n]} \mathbb{E} \left[\|\mathbf{T}_{X_i}\|^3 \right]  < \infty$, 
			\item \label{con:2 Thm} $\lim  \limits_{n \to \infty}\  \max \limits_{i \in [1:n]} \|\mathbf{V}_i\|_{\infty}  < \infty$, and
			\item \label{con:3 Thm} $\mathbf{J} := \lim  \limits_{n \rightarrow \infty} \frac{1}{n} \mathbf{J}_n$ exists and is positive definite with finite determinant, where $ \mathbf{J}_n$ is the FIM in Proposition~\ref{prop:FIM}. 
		\end{enumerate}
		Then, the following two results are true:
		\begin{enumerate}
  \setlength\itemsep{0.01em}
			\item $\boldsymbol{\hat{\theta}}_n$ is a consistent estimator of $\boldsymbol{\theta}_0$, and 
			%
			\item $\sqrt{n} (\boldsymbol{\hat{\theta}}_n-\boldsymbol{\theta}_0) \to \mathcal{N} ( \mathbf{0}_k, \mathbf{J}^{-1})$. 
		\end{enumerate}
	\end{theorem}
  The first two conditions in Theorem~\ref{thm:MainThm} ensure a valid Taylor series expansion of the log-likelihood function and enable one to apply the weak law of large numbers, while the third condition ensures the existence of the asymptotic covariance $\mathbf{J}^{-1}$.
  From Theorem~\ref{thm:MainThm}, it is clear that $\mathbf{J}^{-1}$ is the key to understanding the performance of the MLE.
	\section{Examples}
	\label{sec:Examples}
 In this section, we consider two different densities $p_{X_i}$ that belong to the exponential family. In particular, in Section~\ref{subsec:Example_Gaussian} we focus on the Gaussian distribution, whereas in Section~\ref{subsec:Example_Poisson} we consider the Poisson distribution.
 %
 \subsection{Example~1: Gaussian Distribution} 
 \label{subsec:Example_Gaussian}
 We consider the following model:
\begin{equation}
	X_i = w_i \alpha + C_i, \qquad i \in [1:n],
  \label{eq:systemModel1}
	\end{equation}
where: (i) $\{w_i\}_{i=1}^n$ is a set of known  constants; (ii) $\alpha$ is an unknown deterministic scalar (e.g., the unknown signal in a wireless sensor network~\cite{fang2013}, the amplitude of the radar cross section~\cite{cheng2019target}); and (iii) $C_i \sim \mathcal{N}(0,\sigma^2)$. The model in~\eqref{eq:systemModel1} 
 is widely used in estimation from censored data~\cite{fang2013,khobahi2019,kipnis2022}. 
The probability density function of $X_i$  in~\eqref{eq:systemModel1} is given by
    \begin{align}
        p_{X_i}(x) =  \frac{1}{\sqrt{2 \pi \sigma^2}} \exp\left(-\frac{x^2}{2\sigma^2}+\frac{xw_i\alpha}{\sigma^2} - \frac{w_i^2\alpha^2}{2\sigma^2}\right). 
    \end{align}
We now consider three different cases depending on which quantities we want to estimate.
 
	\noindent {\em{$\bullet$ Case~1: Unknown mean and  known variance.}} With reference to~\eqref{eq:PDFExpFamily}, we have that $\mathrm{T}_x =  x$ and
 \begin{subequations}
\begin{align*}
&h(x) = \frac{1}{\sqrt{2\pi \sigma^2}}\exp \left(-\frac{x^2}{2\sigma^2}\right) , \,\, \phi(\eta_{i,\theta}) =\frac{w_i^2\alpha^2}{2\sigma^2} = \frac{\sigma^2\eta_{i,\theta}^2}{2} , 
\\& \eta_{i,\theta} = v_i\theta \ \text{where} \ v_i = \frac{w_i}{\sigma^2} \ \text{and} \ \theta = \alpha.
\end{align*}
\end{subequations}
For this case, the FIM in \eqref{eq:FIM_censored} is given by (see Appendix~\ref{app:FIM in Case.1})
	\begin{equation}
 \label{eq:FIMGaussCase1}
		 \mathrm{J}_n 
    =  \sum_{i=1}^{n}  w_i^2 \frac{p^2_{X_i}(\tau_i)}{F_{X_i}(\tau_i) (1-F_{X_i}(\tau_i))},
	\end{equation}
where $F_{X_i}$ denotes the cumulative distribution function.
This recovers a well-known result in~\cite{kipnis2022}.
Thus, for Theorem~\ref{thm:MainThm} to hold, it is sufficient that: (i) $\lim_{n\to \infty} \max_{i \in [1:n]} w_i < \infty$, which  often holds in practice; and (ii) $\lim_{n \to \infty} \frac{1}{n}\mathrm{J}_n$ (with $\mathrm{J}_n$ defined in~\eqref{eq:FIMGaussCase1}) exists and is positive, which also holds as long as the $w_i$'s are chosen non-trivially.

The maximum value of the FIM in~\eqref{eq:FIMGaussCase1} occurs at $\tau_i = w_i \alpha$ for $i \in [1:n]$, and it is given by $\mathrm{J}_n = \frac{2}{\pi \sigma^2} \sum_{i=1}^n w_i^2$ (to show this, it suffices to find the first and second derivatives of each term in~\eqref{eq:FIMGaussCase1} with respect to $\tau_i$).
Similarly, the FIM in~\eqref{eq:FIM_uncensored} of estimating $\alpha$ from the uncensored data $\{X_i\}_{i=1}^n$ is given by $\mathrm{I}_n = \frac{1}{\sigma^2} \sum_{i=1}^n w_i^2$. 
Thus, the censoring mechanism requires a fraction $\pi/2$ ($\approx 1.6$) of additional data to match the FIM performance of uncensored data. 
This result extends to sensor fusion, where $n$ is the number of sensors: if $n$ sensors achieve a target FIM with uncensored data, then $\lceil n \pi / 2 \rceil$ sensors are required for the same performance with censored data~\cite{fang2013}.

 
 \noindent {\em{$\bullet$ Case~2: Known mean and unknown variance.}} With reference to~\eqref{eq:PDFExpFamily}, we have that $\mathrm{T}_x \!=\! (x\!-\!\mu_i)^2$ with $\mu_i \!=\! w_i \alpha$, and 
 \begin{align*}
 & \eta_{i,\theta} = v_i\theta \ \text{where} \ v_i = -1/2 \ \text{and} \ \theta = 1/\sigma^2,
 \\
 &h(x) = (2\pi)^{-\frac{1}{2}}, \,\, \,\, \phi(\eta_{i,\theta}) = \log(|\sigma|) = -\frac{1}{2}\log(|2\eta_{i,\theta}|).
 \end{align*}
For this case, the FIM in~\eqref{eq:FIM_censored} is given by (see Appendix~\ref{FIM in Case.2})
	\begin{equation} \label{eq:FIMGaussCase2}
   %
        \mathrm{J}_n = \sum_{i=1}^{n} \frac{\sigma^4}{4} \left(\tau_i - w_i \alpha \right)^2\frac{p_{X_i}^2(\tau_i)}{F_{X_i}(\tau_i)(1-F_{X_i}(\tau_i))}.  
	\end{equation}
 Thus, for Theorem~\ref{thm:MainThm} to hold, similar to Case~1, it suffices that: (i) $\lim_{n\to \infty} \max_{i \in [1:n]} w_i < \infty$, which holds in practice; and (ii) $\lim_{n \to \infty} \frac{1}{n}\mathrm{J}_n$ (with $\mathrm{J}_n$ defined in~\eqref{eq:FIMGaussCase2}) exists and is positive, which, for example, is satisfied  \emph{almost surely} if the $\tau_i$'s are chosen i.i.d.\ from an absolutely continuous distribution.

 \noindent {\em{$\bullet$ Case~3: Unknown mean and unknown variance.}}
	With reference to~\eqref{eq:PDFExpFamily}, we have $\mathbf{T}_x = \begin{bmatrix}
			x & x^2
		\end{bmatrix}^T$, $h(x) = \frac{1}{\sqrt{2\pi }}$, and
\begin{subequations}
\begin{align*}
%
&\boldsymbol{\eta}_{i,\boldsymbol{\theta}}= \mathbf{V}_i \boldsymbol{\theta} \ \text{where} \ \mathbf{V}_i = \begin{bmatrix}
			w_i &0
			\\ 0 & -\frac{1}{2}
		\end{bmatrix} \ \text{and} \ \boldsymbol{\theta} = \begin{bmatrix}
			\frac{\alpha}{\sigma^2} & \frac{1}{\sigma^2}
		\end{bmatrix}^T,
\\& \phi(\boldsymbol{\eta}_{i,\boldsymbol{\theta}}) = \frac{w_i^2\alpha^2}{2\sigma^2} + \log|\sigma|.
\end{align*}
\end{subequations}
For this case, the FIM  in~\eqref{eq:FIM_censored} is given by (see Appendix~\ref{FIM in Case.3})
	\begin{equation}   \label{eq:FIM case3}
	    \begin{aligned}
  %
		\mathbf{J}_n & =\sum_{i=1}^n\sigma^4 \frac{p_{X_i}^2(\tau_i)}{F_{X_i}(\tau_i) (1-F_{X_i}(\tau_i))} 
         \\& \qquad \times
            \begin{bmatrix}
			w_i^2 & -\frac{w_i}{2} \left(\tau_i+w_i\alpha\right)
			\\ -\frac{w_i}{2}\left(\tau_i+w_i\alpha\right)
			&	\frac{(\tau_i+w_i\alpha)^2}{4}
		\end{bmatrix}.
	   \end{aligned}
	\end{equation}
 \begin{remark}
 Each matrix in the sum in~\eqref{eq:FIM case3} is singular. However, this does not necessarily imply that $\mathbf{J}_n$ is singular. To see this, consider the following simple example: $n=2$, $\alpha = 1$, $w_i = 1, i \in [1:2]$, $\sigma=1$, $\tau_1=-1$, and $\tau_2 = 2$. In this case, even if both matrices in the sum in~\eqref{eq:FIM case3} are singular, the determinant of $\mathbf{J}_2$ in~\eqref{eq:FIM case3} is $0.1294$, i.e., $\mathbf{J}_2$ is not singular. 
 \end{remark}
%
Fig.~\ref{Fig: case3 gaussian} offers an illustration of the performance of the MLE for Case~3 when $X_i \sim \mathcal{N}(2,1)$ for all $i \in [1:n]$. For Case~1, the optimal choice of the $\tau_i$'s would be $\tau_i=2$, for all $i \in [1:n]$, whereas for Case~2 the optimal choice of the $\tau_i$'s would be $\tau_i=0.42$ for all $i \in [1:n]$. The dashed curve in Fig.~\ref{Fig: case3 gaussian} corresponds to the case when the $\tau_i$'s in Case~3 are chosen to be either $2$ or $0.42$ with equal probability.
Similarly, the dotted curve in Fig.~\ref{Fig: case3 gaussian} corresponds to the case when the $\tau_i$'s in Case~3 are chosen to be either $1.2$ or $1.9$ with equal probability. These values were selected randomly, but close to the true parameters.
From Fig.~\ref{Fig: case3 gaussian}, we observe that the performance of the MLE improves as $n$ increases; however, for a given $n$, the performance\footnote{Deriving an optimal choice of $\{\tau_i\}_{i=1}^n$ for Case 3 is an interesting open problem, which is worth of further investigation.} depends on the choice of the $\tau_i$'s.
From Fig.~\ref{Fig: case3 gaussian}, it is also apparent that the mean square error from uncensored data (i.e., solid curve) is smaller than the one from censored data.

We conclude the Gaussian example with a proposition providing sufficient conditions on the $\tau_i's$ to satisfy (\emph{almost surely}) condition~3 of Theorem~\ref{thm:MainThm} for Case~3. The proof is in Appendix~\ref{app:Proofprop2}.
\begin{prop}
\label{prop:Case3GaussCond}
Consider the model in~\eqref{eq:systemModel1} with unknown mean and variance. Assume that $\lim_{n\to \infty} \max_{i \in [1:n]} w_i < \infty$,  and  $\lim_{n \to \infty} \frac{| \{w_i: w_i \neq 0 \}|}{n} >0$.
Then, choosing the $\tau_i$'s i.i.d.\ from some absolutely continuous distribution suffices to ensure that $\lim  \limits_{n \rightarrow \infty} \frac{1}{n} \mathbf{J}_n$, with $\mathbf{J}_n$ defined in~\eqref{eq:FIM case3}, is positive definite almost surely.
\end{prop}
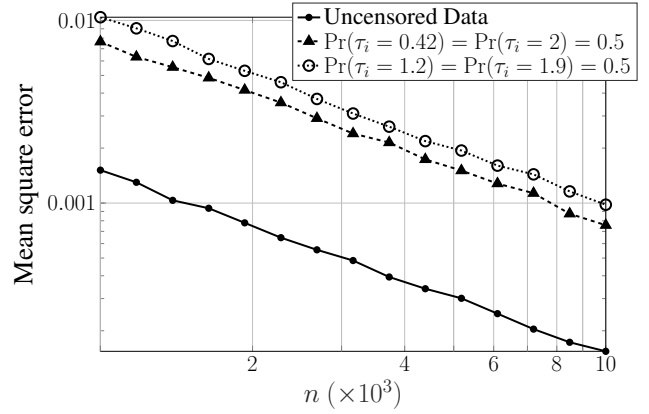
\begin{figure}
    \centering
    \begin{tikzpicture}[scale=.5]
	\begin{axis}[%
		width=5.29in,
		height=3.5in,
		at={(0in,0in)},
scale only axis,
ymode=log,
xmode=log,
xmin=1000,
xmax=10000,
xlabel style={font=\color{white!15!black}\huge},
 xlabel={$n\ (\times 10^{3})$},
ymin=0.000154313343261063,
ymax=0.0104121467485909,
ylabel style={font=\huge},
ylabel={Mean square error},
axis background/.style={fill=white},
xmajorgrids,
ymajorgrids,
minor grid style={opacity=0.1},
		legend style={at={(0.38,0.81)}, anchor=south west, legend cell align=left, align=left, draw=white!15!black,font=\fontsize{18}{20}\selectfont,},
    xticklabel style={font=\LARGE},
    scaled x ticks=false,
    yticklabel style={font=\LARGE},
    log ticks with fixed point,
    xtick={1000,2000,3000,4000,5000,6000,7000,8000,9000,10000},
    xticklabels={,$2$,,$4$,,$6$,,$8$,,{$10$}},
    x tick label style={
           /pgf/number format/1000 sep=
    },
    extra x ticks={1000,3000,5000,7000,9000},
    extra x tick labels={},
    extra x tick style={grid=major, major grid style={opacity=0.3}},
]
\addplot [color=black, line width=1.5pt, mark=*, mark options={solid, black, mark size=2pt}]
  table[row sep=crcr]{%
1000	0.00151448540980237\\
1178	0.00129942426065677\\
1390	0.00103514834066467\\
1638	0.000936592983507007\\
1930	0.000779431806686026\\
2276	0.000646530232037766\\
2682	0.000554690675052597\\
3162	0.000484620737752743\\
3728	0.000393836364814827\\
4394	0.00033968446197376\\
5180	0.000301072508309707\\
6106	0.000248076783746753\\
7196	0.00020415959241114\\
8484	0.00017294137772298\\
10000	0.000154313343261063\\
};
\addlegendentry{Uncensored Data}

\addplot [color=black,dashed,line width=1.5pt, mark=triangle*, mark options={solid, black, mark size=4pt}]
  table[row sep=crcr]{%
1000	0.00764599322408185\\
1178	0.00631334237715115\\
1390	0.00554632224248753\\
1638	0.00484997016698262\\
1930	0.00415754481133309\\
2276	0.00354921533136497\\
2682	0.00291143167286012\\
3162	0.00240163457318798\\
3728	0.00214570546499466\\
4394	0.00173415797064623\\
5180	0.00150638024701537\\
6106	0.00127872900898\\
7196	0.0011324685445038\\
8484	0.000874043340997521\\
10000	0.000756696713115274\\
};
\addlegendentry{$\Pr(\tau_i=0.42) = \Pr(\tau_i=2) = 0.5$}

\addplot [color=black,  dotted, line width=1.5pt, mark=o, mark options={solid, black, mark size=4pt}]
  table[row sep=crcr]{%
1000	0.0104121467485909\\
1178	0.00905111391578067\\
1390	0.00771229746945729\\
1638	0.0061539127299331\\
1930	0.00529760029073589\\
2276	0.00458796877022773\\
2682	0.00371651794888227\\
3162	0.00309259378073108\\
3728	0.0026204936082786\\
4394	0.00218434333512067\\
5180	0.00193958394986315\\
6106	0.0016039407965909\\
7196	0.00143630933870255\\
8484	0.00115837894489054\\
10000	0.000979592247718888\\
};
\addlegendentry{$\Pr(\tau_i=1.2) = \Pr(\tau_i=1.9) = 0.5$}

\end{axis}
\end{tikzpicture}%
    \caption{Case 3: Unknown mean and variance; $X_i \sim \mathcal{N}(2,1)$ for all $i \in [1:n]$.}
    \label{Fig: case3 gaussian}
    \vspace{-0.2cm}
\end{figure}

\subsection{Example~2: Poisson Distribution} 
\label{subsec:Example_Poisson}
We consider the following model: \begin{equation} \label{eq:systemModel2}
X_i \sim \text{Poisson}(\exp(v_i \theta)), \qquad i \in [1:n], 
\end{equation} 
where (i)  $\{v_i\}_{i=1}^n$ is a set of known finite constants, and (ii) $\theta$ is an unknown deterministic scalar. 
This model is widely used in low-power laser communications~\cite{dytsoPoisson}, where $\exp(v_i\theta)$ relates to the applied current, and in photonics applications like image acquisition with binary Poisson statistics~\cite{yang2011bits}, where $v_i$ represents the detector efficiency and $\theta$ is the light source intensity.
The probability mass function of $X_i$ belongs to the exponential family and the FIM in~\eqref{eq:FIM_censored} is given by
\begin{equation}
 \label{eq:FIMPoisson}
	\mathrm{J}_n \!=\!  \sum_{i=1}^{n}  v_i^2 
            \exp (2v_i\theta) \frac{\left(F_{X_i}( \tau_i ) \!-\! F_{X_i}( \tau_i -1)\right)^2}{F_{X_i}( \tau_i )\left(1-F_{X_i}( \tau_i )\right)} ,
	\end{equation}
which is computed in Appendix~\ref{app:FIMpoisson}. For Theorem~\ref{thm:MainThm} to hold, it suffices that: (i) $\lim_{n\to \infty} \max_{i \in [1:n]} v_i < \infty$, which holds in practice; 
and (ii) $\lim_{n \to \infty} \frac{1}{n}\mathrm{J}_n$ (with $\mathrm{J}_n$ defined in~\eqref{eq:FIMPoisson}) exists and is positive, which also holds as long as the $\tau_i \geq 0$ for all $i \in [1:n]$ and the $v_i$'s are chosen non-trivially.

	
	\section{CONCLUSION} \label{sec:conclusion}
This work has established regularity conditions for the \emph{consistency} and \emph{asymptotic normality} of the MLE of multiple parameters from a censored GLM model. The derived conditions can be considered mild. 
 The FIM  under both censored and uncensored data scenarios has been characterized, which is crucial for quantifying the information loss due to censoring. Finally, the asymptotic performance of the MLE has been analyzed for two common noise distributions: Gaussian, with both unknown mean and variance, and Poisson, with an unknown mean.

    \appendices
	\section{Proof of~(5)} 
 \label{app:ProofSimplifLogLik}
	From the definition of the log-likelihood function, 
	\begin{align} 
		&\ell_n \left (\boldsymbol{\theta};\{b_i\}_{i=1}^n \right ) 
		%
		= \sum_{i=1}^n \log \left( P_{B_i}(b_i;\boldsymbol{\theta}) \right ) \label{eq:loglikelihood}
		\\& = \sum_{i=1}^n \left [ \log \left( P_{B_i}(b_i;\boldsymbol{\theta}) \right ) \!+\! \phi(\boldsymbol{\eta}_{i,\boldsymbol{\theta}}) \!-\! \phi(\boldsymbol{\eta}_{i,\boldsymbol{\theta}}) \right ] \notag
		\\& \stackrel{{\rm{(a)}}}{=} \sum_{i=1}^n \Big[ \log \Big( \int \limits_{\mathcal{X}(b_i)} p_{X_i}(x;\boldsymbol{\theta} ) \ {\rm{d}}x \Big ) + \phi(\boldsymbol{\eta}_{i,\boldsymbol{\theta}}) - \phi(\boldsymbol{\eta}_{i,\boldsymbol{\theta}}) \Big ] \notag
		\\& \stackrel{{\rm{(b)}}}{=} \sum_{i=1}^n \Big[ \log \Big( \int \limits_{\mathcal{X}(b_i)} h(x) \exp (\langle \boldsymbol{\eta}_{i,\boldsymbol{\theta}}, \mathbf{T}_x \rangle - \phi(\boldsymbol{\eta}_{i,\boldsymbol{\theta}})) \ {\rm{d}}x \Big )  \notag
		\\& \qquad \qquad   + \phi(\boldsymbol{\eta}_{i,\boldsymbol{\theta}}) - \phi(\boldsymbol{\eta}_{i,\boldsymbol{\theta}}) \Big ] \notag 
		\\& = \sum_{i=1}^n \Big[ \log \Big( \int \limits_{\mathcal{X}(b_i)} h(x) \exp (\langle \boldsymbol{\eta}_{i,\boldsymbol{\theta}}, \mathbf{T}_x \rangle ) \ {\rm{d}}x \Big )  \!-\! \phi(\boldsymbol{\eta}_{i,\boldsymbol{\theta}}) \Big] \notag 
		\\& \stackrel{{\rm{(c)}}}{=}   \sum \limits_{i=1}^n \left[	\phi(\boldsymbol{\eta}_{i,\boldsymbol{\theta}};b_i) - \phi(\boldsymbol{\eta}_{i,\boldsymbol{\theta}})\right], \label{eq:ll_final}
	\end{align}
	where the labeled equalities follow from:
	$\rm{(a)}$ substituting~\eqref{eq:PMFBi};
	$\rm{(b)}$ using~\eqref{eq:PDFExpFamily};
	and $\rm{(c)}$ letting $\phi(\boldsymbol{\eta}_{i,\boldsymbol{\theta}};b_i)$ be the log-partition function of $p_{X_i | B_i}(x|b_i)$ (note that $\phi(\boldsymbol{\eta}_{i,\boldsymbol{\theta}};b_i)$ is a normalization quantity that ensures that $p_{X_i | B_i}(x|b_i)$ is a valid density).
	
	\section{Proof of Proposition~1}
	\label{app:FIM}
	For brevity, let $\theta_r$ denote the $r^{th}$ element of $\boldsymbol{\theta}$ for all $r \in [1:k]$. 
	The $(r,s)^{th}$ element of the FIM, indicated as $\mathbf{J}_{r,s}$, is
	\begin{align}
		&\mathbf{J}_{r,s} = \mathbb{E}\Big[ \sum_{i=1}^n  \sum_{j=1}^n\frac{\partial \log P_{B_i}(B_i; \boldsymbol{\theta})}{\partial \theta_{r}} \frac{\partial \log P_{B_j}(B_j; \boldsymbol{\theta})}{\partial \theta_{s}}\Big] \notag \\
		&\stackrel{{\rm{(a)}}}{=} \mathbb{E}\Big[ \sum_{i=1}^n  \frac{\partial \log P_{B_i}(B_i; \boldsymbol{\theta})}{\partial \theta_{r}} \frac{\partial \log P_{B_i}(B_i; \boldsymbol{\theta})}{\partial \theta_{s}}\Big] \notag \\
		& \stackrel{{\rm{(b)}}}{=} \sum_{i=1}^{n} \mathbb{E} \Big [\Big( \sum_{j=1}^d \left(\mathbb{E} \left [ \mathbf{T}_{X_i}(j)|B_i\right ] - \mathbb{E}\left[\mathbf{T}_{X_i}(j) \right ] \right ) \mathbf{V}_i(j,r) \Big )  \notag
		\\& \qquad \qquad  \Big( \sum_{\ell=1}^d \left(\mathbb{E} \left [ \mathbf{T}_{X_i}(\ell)|B_i\right ] - \mathbb{E}\left[\mathbf{T}_{X_i}(\ell) \right ] \right ) \mathbf{V}_i(\ell,s)\Big ) \Big ] \notag \\
		& = \sum_{i=1}^{n} \sum_{j=1}^d \sum_{\ell=1}^d \mathsf{Cov} \left( \mathbb{E} \left [ \mathbf{T}_{X_i}(j)|B_i\right ],\mathbb{E} \left [ \mathbf{T}_{X_i}(\ell)|B_i\right ] \right ) \notag
		\\& \hspace{4cm}  \times\mathbf{V}_i(j,r)\mathbf{V}_i(\ell,s) \notag \\
		&\stackrel{{\rm{(c)}}}{=} \sum_{i=1}^n \mathbf{V}_i^T(r) \mathsf{Cov} \left( \mathbb{E} \left [ \mathbf{T}_{X_i}|B_i\right ] \right ) \mathbf{V}_i(s),
	\end{align}
	where the labeled equalities follow from:
	${\rm{(a)}}$ the fact that the $B_i$'s are independent; hence, when $i \neq j$, the expected value of the product can be written as the product of the expected values and the expected value of the score is zero (see Lemma~\ref{lem:Score} below);
	$\rm{(b)}$ applying Lemma~\ref{lem:Score} below;
 and $\rm{(c)}$ letting $\mathbf{V}_i(s)$ denote the $s$-th column of $\mathbf{V}_i$.
	This concludes the proof of Proposition~\ref{prop:FIM}.
	
	\begin{lem}
		\label{lem:Score}
		For every $i \in [1:n]$ and $r \in [1:k]$, the score function is given by
		\begin{equation}
			\frac{\partial \log P_{B_i}}{\partial \theta_{r}} \!=\!\! \sum_{j=1}^d \!\left(\mathbb{E} \!\left [ \mathbf{T}_{X_i}(j)|B_i\right ] \!-\! \mathbb{E}\left[\mathbf{T}_{X_i}(j) \right ] \right ) \!\mathbf{V}_i(j,r),
		\end{equation}
  where $\mathbf{T}_{X_i}(j)$ is the $j$-th element of $\mathbf{T}_{X_i}$ and $\mathbf{V}_i(j,r)$ is the $(j,r)$-th element of $\mathbf{V}_i$.
	\end{lem}
	\begin{proof}
		From~\eqref{eq:ll_final} we have $\log \left( P_{B_i}(b_i;\boldsymbol{\theta}) \right ) = \phi(\boldsymbol{\eta}_{i,\boldsymbol{\theta}};b_i) - \phi(\boldsymbol{\eta}_{i,\boldsymbol{\theta}})$ and by applying the chain rule, we have that
		\begin{align}
			& \frac{\partial \log P_{B_i}(b;\boldsymbol{\theta})}{\partial \theta_r} = 
			\frac{\partial}{\partial \theta_r} \left( \phi(\boldsymbol{\eta}_{i,\boldsymbol{\theta}};b) - \phi(\boldsymbol{\eta}_{i,\boldsymbol{\theta}}) \right )  
   \notag
			\\& = \sum_{j=1}^d \left( \frac{\partial \phi(\boldsymbol{\eta}_{i,\boldsymbol{\theta}};b)}{\partial \boldsymbol{\eta}_{i,\boldsymbol{\theta}} (j)} - \frac{\partial \phi(\boldsymbol{\eta}_{i,\boldsymbol{\theta}})}{\partial \boldsymbol{\eta}_{i,\boldsymbol{\theta}} (j)} \right ) \frac{\partial \boldsymbol{\eta}_{i,\boldsymbol{\theta}} (j)}{\partial \theta_r} \notag
			\\& = \sum_{j=1}^d \left( \frac{\partial \phi(\boldsymbol{\eta}_{i,\boldsymbol{\theta}};b)}{\partial \boldsymbol{\eta}_{i,\boldsymbol{\theta}} (j)} - \frac{\partial \phi(\boldsymbol{\eta}_{i,\boldsymbol{\theta}})}{\partial \boldsymbol{\eta}_{i,\boldsymbol{\theta}} (j)} \right ) \mathbf{V}_i(j,r) \notag
			\\& = \sum_{j=1}^d \left(\mathbb{E} \left [ \mathbf{T}_{X_i}(j)|B_i = b\right ] - \mathbb{E}\left[\mathbf{T}_{X_i}(j) \right ] \right ) \mathbf{V}_i(j,r),
			\label{eq:FirstDerivativeMultivariate}
		\end{align}
		where the last equality follows from~\cite{dasgupta2011exponential}. 
  Note also that by the law of total expectation, we have $\mathbb{E}\left [ \frac{\partial \log P_{B_i}}{\partial \theta_{r}}\right ] =0$.
  This concludes the proof of Lemma~\ref{lem:Score}.
	\end{proof}
	
	\section{Proof of Theorem~1}
	\label{app:MainTheorem}
In~\cite{bradley1962}, the authors developed a theory for the MLE when the observations are independent and come from distinct, yet related populations, i.e., with some parameters in common. The authors referred to such populations as {\em {associated}}. In particular, the authors derived
 regularity conditions under which the MLE of parameters in associated populations is shown to be consistent and asymptotically normal.
 In what follows, we tailor these conditions to our GLM with 1-bit measurements described in Section~\ref{sec:ProbForm}.

	\begin{enumerate}[label={\arabic*)}]
		\item \label{con:1.1_multi}
		We require the existence of the following partial derivatives to ensure that the Taylor series expansion of the log-likelihood function in~\eqref{eq:loglikelihood} exists~\cite[conditions I(i)]{bradley1962}, 
		\begin{equation}
			\frac{\partial \log P_{B_i}}{\partial \theta_r}, \ \frac{\partial^2 \log P_{B_i}}{\partial \theta_r \partial \theta_s}, \ \text { and } \ \frac{\partial^3 \log P_{B_i}}{\partial \theta_r \partial \theta_s \partial \theta_t}
		\end{equation}
		for all $(r, s, t) \in [1:k]^3$ and $i \in [1:n]$.
		We start with the first order partial derivative. 
		From Lemma~\ref{lem:Score},
		\begin{align*}
			&\frac{\partial \log P_{B_i}}{\partial \theta_{r}} 
			 \!=\!\! \sum_{j=1}^d \!\left(\mathbb{E} \left [ \mathbf{T}_{X_i}(j)|B_i\right ] \!-\! \mathbb{E}\left[\mathbf{T}_{X_i}(j) \right ] \right ) \!\mathbf{V}_i(j,r).
		\end{align*}
		Thus, condition~1 and condition~2 in Theorem~\ref{thm:MainThm} ensure the existence of $ \partial \log P_{B_i} / \partial \theta_r$ for all $r \in [1:k]$.
		
		Similarly, for the second order partial derivative, by applying the chain rule, we arrive at
		\begin{align}
			&\frac{\partial^2 \log P_{B_i}}{\partial \theta_r \partial \theta_s} = 
			\frac{\partial^2}{\partial \theta_r \partial \theta_s} \left( \phi(\boldsymbol{\eta}_{i,\boldsymbol{\theta}};b) - \phi(\boldsymbol{\eta}_{i,\boldsymbol{\theta}}) \right ) 
			\notag
			\\& = \sum_{j=1}^d \sum_{\ell=1}^d  \left( \frac{\partial^2 \phi(\boldsymbol{\eta}_{i,\boldsymbol{\theta}};b)}{\partial \boldsymbol{\eta}_{i,\boldsymbol{\theta}} (\ell) \partial \boldsymbol{\eta}_{i,\boldsymbol{\theta}} (j)} - \frac{\partial^2 \phi(\boldsymbol{\eta}_{i,\boldsymbol{\theta}})}{\partial \boldsymbol{\eta}_{i,\boldsymbol{\theta}} (\ell)\partial \boldsymbol{\eta}_{i,\boldsymbol{\theta}} (j)} \right ) \notag
			\\ & \qquad \qquad \times \mathbf{V}_i(\ell,s) \mathbf{V}_i(j,r)  \notag
			\\& = \sum_{j=1}^d \sum_{\ell=1}^d \left[ \mathsf{Cov} \left( \mathbf{T}_{X_i}(j),\mathbf{T}_{X_i}(\ell) | B_i = b\right ) \right. \notag
			\\
			& \left. \quad - \mathsf{Cov} \left( \mathbf{T}_{X_i}(j),\mathbf{T}_{X_i}(\ell) \right ) \right ] \mathbf{V}_i(\ell,s) \mathbf{V}_i(j,r), \label{eq:SecondDerivativeMulti}
		\end{align}
		where the last equality follows from~\cite{dasgupta2011exponential}.
  Thus, condition 1 and condition 2 in Theorem~\ref{thm:MainThm} ensure the existence of $\partial^2 \log P_{B_i}/\partial{\theta_r}\partial{\theta_s}$ for all $(r,s)\in [1:k]^2$.
		
		Finally, for the third order partial derivative, by applying the chain rule, we have that
		\begin{align}
			& \frac{\partial^3 \log P_{B_i}}{\partial \theta_r \partial \theta_s \partial \theta_t} = 
			\frac{\partial^3}{\partial \theta_r \partial \theta_s \partial \theta_t} \left( \phi(\boldsymbol{\eta}_{i,\boldsymbol{\theta}};b) - \phi(\boldsymbol{\eta}_{i,\boldsymbol{\theta}}) \right ) 
			\notag
			\\& = \sum_{j=1}^d \sum_{\ell=1}^d \sum_{k=1}^d  \left[ \frac{\partial^3 \phi(\boldsymbol{\eta}_{i,\boldsymbol{\theta}};b)}{\partial \boldsymbol{\eta}_{i,\boldsymbol{\theta}} (k) \partial \boldsymbol{\eta}_{i,\boldsymbol{\theta}} (\ell) \partial \boldsymbol{\eta}_{i,\boldsymbol{\theta}} (j)} \right. \notag
			\\ & \left.- \frac{\partial^3 \phi(\boldsymbol{\eta}_{i,\boldsymbol{\theta}})}{\partial \boldsymbol{\eta}_{i,\boldsymbol{\theta}} (k) \partial \boldsymbol{\eta}_{i,\boldsymbol{\theta}} (\ell)\partial \boldsymbol{\eta}_{i,\boldsymbol{\theta}} (j)}\! \right ]\! \mathbf{V}_i(k,t) \! \mathbf{V}_i(\ell,s) \! \mathbf{V}_i(j,r)  \notag
			\\& = \sum_{j=1}^d \sum_{\ell=1}^d \sum_{k=1}^d   \left[ \kappa_{\mathbf{T}_{X_i|B_i=b}}(j, \ell,k)  - \kappa_{\mathbf{T}_{X_i}}(j, \ell,k) \right] \notag
			\\& \qquad \qquad \times  \mathbf{V}_i(k,t) \mathbf{V}_i(\ell,s) \mathbf{V}_i(j,r),
			\label{eq:third_derivative_final}
		\end{align}
		where the last equality follows from~\cite{dasgupta2011exponential} where
		\begin{subequations}
			\begin{align}
				&\kappa_{\mathbf{T}_{X_i|B_i=b}}(j, \ell,k) = \notag\\
    &\mathbb{E}\left[\prod_{u \in \{j,\ell,k\}}(\mathbf{T}_{X_i}(u)  -  \mathbb{E}[\mathbf{T}_{X_i}(u)|B_i=b])|B_i=b  \right],
				\\& \kappa_{\mathbf{T}_{X_i}}(j, \ell,k) = \mathbb{E}\left[\prod_{u \in \{j,\ell,k\}}(\mathbf{T}_{X_i}(u) - \mathbb{E}[\mathbf{T}_{X_i}(u)] ) \right].
			\end{align}
		\end{subequations}
		Thus, condition 1 and condition 2 in Theorem~\ref{thm:MainThm} ensure the existence of $\partial^3 \log P_{B_i}/\partial{\theta_r} \partial{\theta_s}\partial{\theta_t}$ for all $(r,s,t)\in [1:k]^3$.
		\item \label{con:1.2 multi}
		Now we check~\cite[conditions I(ii)]{bradley1962}, which consist of two parts.
  First, we require the convergence of the first and second order partial derivatives, which will allow to interchange the differentiation and summation. This is indeed satisfied, as from~\eqref{eq:BinaryRV} we have
		\begin{align}
			& \sum_{b \in \{1,-1\}} \frac{\partial P_{B_i}(b;\boldsymbol{\theta})}{\partial \theta_r} = \frac{\partial P_{B_i}(1;\boldsymbol{\theta})}{\partial \theta_r} + \frac{\partial P_{B_i}(-1;\boldsymbol{\theta})}{\partial \theta_r}  \notag
			\\ & = \frac{\partial F_{X_i}(\tau_i)}{\partial \theta_r} + \frac{\partial (1 - F_{X_i}(\tau_i))}{\partial \theta_r} = 0,
		\end{align}
		and similarly we have that
		\begin{equation}
			\sum_{b \in \{1,-1\}}\hspace{-.3cm} \frac{\partial^{2} P_{B_i}\!(b;\boldsymbol{\theta})}{\partial \theta_{r} \partial \theta_{s}} 
			\!=\!\frac{\partial^2 F_{X_i}\!(\tau_i)}{\partial \theta_{r} \partial \theta_{s}}\!+\! \frac{\partial^2(1\!-\! F_{X_i}(\tau_i))}{\partial \theta_{r} \partial \theta_{s}}\!=\!0.
		\end{equation}
		Hence, both the first order and second order partial derivatives converge.
		
		Second, we need to check that the third order derivative is finite. From~\eqref{eq:third_derivative_final}, we have that
		\begin{align}
			& \frac{\partial^3 \log P_{B_i}}{\partial \theta_r \partial \theta_s \partial \theta_t} = \sum_{j=1}^d \sum_{\ell=1}^d \sum_{k=1}^d   \Big[ \kappa_{\mathbf{T}_{X_i|B_i=b}}(j, \ell,k)  \notag
			\\ & - \kappa_{\mathbf{T}_{X_i}}(j, \ell,k) \Big] \mathbf{V}_i(k,t) \mathbf{V}_i(\ell,s) \mathbf{V}_i(j,r). 
		\end{align}
		Under condition~\ref{con:1 Thm} and condition~\ref{con:2 Thm} in Theorem~\ref{thm:MainThm}, there exists some finite positive constant $K$ such that
		\begin{align}
			& \Big | \sum_{j=1}^d \sum_{\ell=1}^d \sum_{k=1}^d   \left[ \kappa_{\mathbf{T}_{X_i|B_i=b}}(j, \ell,k)  - \kappa_{\mathbf{T}_{X_i}}(j, \ell,k) \right] \notag
			\\ & \qquad \times \mathbf{V}_i(k,t) \mathbf{V}_i(\ell,s) \mathbf{V}_i(j,r) \Big | \leq K,
		\end{align}
		for all $(r, s, t) \in [1:k]^3$ and $i \in [1:n]$.
		In other words, with reference to~\cite[condition I(ii)]{bradley1962}, we have that $H_{i r s t}(b) = K$ for all $(r, s, t) \in [1:k]^3$, $b \in \{-1,1\}$, and $i \in [1:n]$.
		Now, it holds that
		\begin{equation}
			\sum_{b \in\{1,-1\}}\hspace{-.35cm}H_{i r s t}\!\left(b\right)\! P_{B_i}\!(b;\boldsymbol{\theta})\!=\! K\!\!\!\!  \sum_{b \in\{1,-1\}}\!\!\!\! P_{B_i}(b;\boldsymbol{\theta})\! =\! K,
		\end{equation}
		for all $\boldsymbol{\theta} \in \Theta$, $(r, s, t) \in [1:k]^3$, and $i \in [1:n]$. Hence, for all $i \in [1:n]$, with reference to~\cite[condition I(ii)]{bradley1962}, we can set $M_i=K$.
		\item \label{con:2.1_multi}
This condition is~\cite[condition II(i)]{bradley1962} and it is needed to apply the weak law of large numbers for independent random variables~\cite[p.174]{cramer1947problems}.  
		In order to check this condition we first need to find the set $D_{1 i}$, which for all $i \in [1:n]$ and for all $r \in [1:k]$ is given by
		\begin{equation}
			D_{1 i} \!=\! \left\{b \in \{-1,1\}: \left| \frac{\partial \log P_{B_i}(b; \boldsymbol{\theta})}{ \partial \theta_r} \right|>n  \right\}. 	
		\end{equation}
  Under condition~\ref{con:1 Thm} and condition~\ref{con:2 Thm} in Theorem~\ref{thm:MainThm}, the derivative $ \partial \log P_{B_i} / \partial \theta_r$ is finite for all $ i \in [1:n]$ and $r \in [1:k]$. Thus, for a sufficiently large $n$, there will be no $b \in\{-1,1\}$ such that $\left|\frac{\partial \log P_{B_i}(b;\boldsymbol{\theta})}{\partial \theta}\right|>n$.
		Hence, for large $n$, the set $D_{1 i}$ will be empty, i.e.,
		\begin{equation}
			D_{1 i}= \varnothing \quad \text{ for all } i \in[1:n].
		\end{equation}
		Hence, 
		\begin{equation}
			\sum_{i=1}^n \sum_{b \in D_{1 i}} P_{B_i}\left(b;\boldsymbol{\theta}\right)=0,
		\end{equation}
		which satisfies the condition $\sum_{i=1}^n \sum_{b \in D_{1 i}} P_{B_i}\left(b;\boldsymbol{\theta}\right)=o(1)$~\cite[conditions II(i)]{bradley1962}.
		Similarly, we can find the set $D_{2 i}$ for all $i \in [1:n]$ and $r \in [1:k]$ as follows,
		\begin{equation}
			D_{2 i} \!=\! \left\{b \in \{-1,1\}: \left| \frac{\partial \log P_{B_i}(b;\boldsymbol{\theta})}{ \partial \theta_r} \right| < n\right\} . 
		\end{equation}
		Using the same reasoning as above, $D_{2i} = \{-1,1\}$ for a sufficiently large $n$ whenever condition~\ref{con:1 Thm} and condition~\ref{con:2 Thm} in Theorem~\ref{thm:MainThm} hold.
		Now, we need to check the following condition~\cite[condition II(i)]{bradley1962},
		\begin{equation}
			\sum_{i=1}^{n} \sum_{b \in D_{2 i}}\!\left(\frac{\partial \log P_{B_i}(b; \boldsymbol{\theta})}{\partial \theta_r}\!\right)^{2} \!\!P_{B_i}(b;\boldsymbol{\theta}) = o(n^{2}).
   \label{eq:IntermStep3}
		\end{equation}
		From Lemma~\ref{lem:Score}, we have that
		\begin{equation}
			\frac{\partial\log P_{B_i}}{\partial \theta_r} \!=\!\sum_{j=1}^d\!\left(\mathbb{E}[ \mathbf{T}_{X_i}\!(j)|B_i ]- \mathbb{E}[\mathbf{T}_{X_i}(j) ] \right)\!\mathbf{V}_i(j,r).
		\end{equation}
		Hence, whenever condition~\ref{con:1 Thm} and condition~\ref{con:2 Thm} in Theorem~\ref{thm:MainThm} hold, there exists some constant $C$ such that
		\begin{equation}
			\sum_{j=1}^d \!\left(\mathbb{E} [ \mathbf{T}_{X_i}(j)|B_i ]\!-\! \mathbb{E}[\mathbf{T}_{X_i}(j)  ] \right ) \mathbf{V}_i(j,r)\leq \sqrt{C},
		\end{equation}
		for $b \in \{-1,1\}$ and for all $i \in [1:n]$ and $r \in [1:k]$.
		Therefore, the left-hand side of~\eqref{eq:IntermStep3} can be written as
		\begin{align}
			& \sum_{i=1}^n \sum_{b \in\{-1,1\}} \Big( \sum_{j=1}^d \left(\mathbb{E} \left [ \mathbf{T}_{X_i}(j)|B_i=b\right ] \right.  \notag
			\\ & \qquad \qquad  \left.
			- \mathbb{E}\left[\mathbf{T}_{X_i}(j) \right ] \right ) \mathbf{V}_i(j,r)\Big )^2  P_{B_i}\left(b;\boldsymbol{\theta}\right)  \notag
			\\& \leq \sum_{i=1}^n \sum_{b \in\{-1,1\}} C P_{B_i}\left(b;\boldsymbol{\theta}\right) 
			=C n.
		\end{align}
		Thus,~\eqref{eq:IntermStep3} is satisfied for all $r \in [1:k]$.
		\item \label{con:2.2_multi}
  This condition is~\cite[condition II(ii)]{bradley1962} and, as the one above, it is needed to apply the weak law of large numbers for independent random variables~\cite[p.174]{cramer1947problems}.
		In order to check this condition we first need to find the set $D_{3 i}$  which for all $i \in [1:n]$, is given by		
		\begin{equation}
			D_{3 i} \!=\! \left\{b \!\in\! \{-1,1\}: \left|\frac{\partial^{2} \log P_{B_i}(b;\boldsymbol{\theta})}{ \partial \theta_r \partial \theta_s}\right|>n  \right\}.
		\end{equation}
		Whenever condition~\ref{con:1 Thm} and condition~\ref{con:2 Thm} in Theorem~\ref{thm:MainThm} hold, the derivative $ \partial^2 \log P_{B_i} / \partial \theta_r \partial \theta_s$ is finite for all $i \in [1:n]$ (see our analysis of~\eqref{eq:SecondDerivativeMulti}). 
		Thus, for a sufficiently large $n$, there will be no $b \in\{-1,1\}$ such that $\left|\frac{\partial^2 \log P_{B_i}(b;\boldsymbol{\theta})}{\partial \theta_r \partial \theta_s}\right|>n$. Hence, for large $n$, the set $D_{3i}$ will be empty, i.e., 
		\begin{equation}
			D_{3i}= \varnothing \quad \text{ for all } \quad  i \in[1:n].
		\end{equation}
		Therefore, 
		\begin{equation}
			\sum_{i=1}^{n} \sum_{b \in D_{3 i}} P_{B_i}(b;\boldsymbol{\theta})=0 = o(1),
		\end{equation}
  as required by~\cite[condition II(ii)]{bradley1962}.
		Similarly, we can find the set $D_{4 i}$ for all $i \in [1:n]$ as follows,
		\begin{equation}
			D_{4 i}  \!=\! \left\{b \!\in\! \{-1,1\}: \left| \frac{\partial^2 \log P_{B_i}(b;\boldsymbol{\theta})}{ \partial \theta_r \partial \theta_r} \right| < n \right\}. 
		\end{equation}
		Using the same reasoning as before we have, $D_{4 i} = \{-1,1\}$ for a sufficiently large $n$.
		Now, we need to check the following condition~\cite[condition II(ii)]{bradley1962},
		\begin{equation}
  \label{eq:IntermStep4}
			\sum_{i=1}^{n} \!\sum_{b \in D_{4 i}}\!\!\left(\!\frac{\partial^{2} \log P_{B_i}(b;\boldsymbol{\theta})}{\partial \theta_{r} \partial \theta_{s}}\!\right)^{2} \!\! P_{B_i}(b;\boldsymbol{\theta}) \!=\! o(n^{2}),  
		\end{equation}
  where $\frac{\partial^2 \log P_{B_i}}{\partial \theta_r \partial \theta_s}$ is given in~\eqref{eq:SecondDerivativeMulti}.
		Hence, whenever condition~\ref{con:1 Thm} and condition~\ref{con:2 Thm} in Theorem~\ref{thm:MainThm} hold, there exists some constant $\tilde{C}$ such that
		\begin{align}
			& \sum_{j=1}^d \sum_{\ell=1}^d \Big( \mathsf{Cov} \left( \mathbf{T}_{X_i}(j),\mathbf{T}_{X_i}(\ell) | B_i = b\right )  \notag
			\\ & - \mathsf{Cov} \left( \mathbf{T}_{X_i}(j),\mathbf{T}_{X_i}(\ell) \right ) \Big ) \mathbf{V}_i(\ell,s) \mathbf{V}_i(j,r) \leq \sqrt{\tilde{C}},
		\end{align}
		for all $i \in [1:n]$ and $(r,s)\in [1:k]^2$.
		Therefore, 
		\begin{align}
			& \sum_{i=1}^{n} \sum_{b \in D_{4 i}}\left(\frac{\partial^{2} \log P_{B_i}(b;\boldsymbol{\theta})}{\partial \theta_r \partial \theta_s}\right)^{2} P_{B_i}(b;\boldsymbol{\theta}) \notag
			\\ & \leq \sum_{i=1}^n \sum_{b \in\{-1,1\}} \tilde{C} P_{B_i}\left(b; \boldsymbol{\theta}\right) = 
			\tilde{C}  n.
		\end{align}
\begin{figure*}[h]
			\begin{align}
				&\lim _{n \rightarrow \infty} \frac{1}{n} \sum_{i=1}^{n} \sum_{b \in \{1,-1\} } \frac{\partial \log P_{B_i}(b; \boldsymbol{\theta})}{\partial \theta_{r}} \frac{\partial \log P_{B_i}(b; \boldsymbol{\theta})}{\partial \theta_{s}} P_{B_i}(b; \boldsymbol{\theta})  =\lim _{n \rightarrow \infty} \frac{1}{n} \sum_{i=1}^{n} \mathbb{E} \left [\frac{\partial \log P_{B_i}(B_i; \boldsymbol{\theta})}{\partial \theta_{r}} \frac{\partial \log P_{B_i}(B_i; \boldsymbol{\theta})}{\partial \theta_{s}} \right ] \notag
				\\& =\lim _{n \rightarrow \infty} \frac{1}{n} \sum_{i=1}^{n} \mathbb{E} \left [\left( \sum_{j=1}^d \left(\mathbb{E} \left [ \mathbf{T}_{X_i}(j)|B_i\right ] - \mathbb{E}\left[\mathbf{T}_{X_i}(j) \right ] \right ) \mathbf{V}_i(j,r) \right ) \left( \sum_{j=1}^d \left(\mathbb{E} \left [ \mathbf{T}_{X_i}(j)|B_i\right ] - \mathbb{E}\left[\mathbf{T}_{X_i}(j) \right ] \right ) \mathbf{V}_i(j,s)\right ) \right ] \notag
				\\& = \lim _{n \rightarrow \infty} \frac{1}{n} \sum_{i=1}^{n} \sum_{j=1}^d \sum_{\ell=1}^d \mathsf{Cov} \left( \mathbb{E} \left [ \mathbf{T}_{X_i}(j)|B_i\right ],\mathbb{E} \left [ \mathbf{T}_{X_i}(\ell)|B_i\right ] \right )\mathbf{V}_i(j,r)\mathbf{V}_i(\ell,s)
				\label{eq:cond 2.3 RHSMulti}
			\end{align}
   \hrule
   \begin{align}
				& \lim _{n \rightarrow \infty} \frac{1}{n} 	\sum_{i=1}^{n} \sum_{b \in D_{4 i}}-\left(\frac{\partial^{2} \log P_{B_i}(b; \boldsymbol{\theta})}{\partial \theta_{r} \partial \theta_{s}}\right) P_{B_i}(b;\boldsymbol{\theta}) \stackrel{{\rm{(a)}}}{=} \lim _{n 	\rightarrow \infty} \frac{1}{n} \sum_{i=1}^{n} - \mathbb{E} \left [\frac{\partial^{2} \log P_{B_i}(b; \boldsymbol{\theta})}{\partial \theta_{r} \partial \theta_{s}} \right ] \notag
				\\& = \lim _{n \rightarrow \infty} 	\frac{1}{n} \sum_{i=1}^{n}  \sum_{j=1}^d \sum_{\ell=1}^d \left(\mathsf{Cov} \left( \mathbf{T}_{X_i}(j),\mathbf{T}_{X_i}(\ell) \right ) - \mathbb{E} \left [ \mathsf{Cov} \left( \mathbf{T}_{X_i}(j),\mathbf{T}_{X_i}(\ell) | B_i = b \right )\right ] \right ) \mathbf{V}_i(\ell,s) \mathbf{V}_i(j,r) \notag
				\\& \stackrel{{\rm{(b)}}}{=} \lim _{n 	\rightarrow \infty} \frac{1}{n} \sum_{i=1}^{n} \sum_{j=1}^d \sum_{\ell=1}^d \mathsf{Cov} \left( \mathbb{E} \left [ \mathbf{T}_{X_i}(j)|B_i\right ],\mathbb{E} \left [ \mathbf{T}_{X_i}(\ell)|B_i\right ] \right )\mathbf{V}_i(j,r)\mathbf{V}_i(\ell,s)
				\label{eq:cond 2.3 LHSMulti}
			\end{align}
   \hrule
\end{figure*} 

Thus,~\eqref{eq:IntermStep4}
		is satisfied.
  Always with reference to~\cite[condition II(ii)]{bradley1962}, we also need to verify that the two limits in~\cite[eq.13]{bradley1962} are the same.
  In particular, the limit on the right-hand side of~\cite[eq.13]{bradley1962} is given in~\eqref{eq:cond 2.3 RHSMulti}, whereas the limit on the left-hand side of~\cite[eq.13]{bradley1962} is given in~\eqref{eq:cond 2.3 LHSMulti}, both at the top of the next page, where the equality in $\rm{(a)}$ follows since, as proved above, for a sufficiently large $n$, we have that $D_{4 i} = \{-1,1\}$; and the equality in $\rm{(b)}$ is due to the law of total covariance.

		


		
		From~\eqref{eq:cond 2.3 RHSMulti} and~\eqref{eq:cond 2.3 LHSMulti}, it follows that the two limits in~\cite[eq.13]{bradley1962} are the same.
		Moreover, from~\eqref{eq:cond 2.3 LHSMulti} (or equivalently~\eqref{eq:cond 2.3 RHSMulti}), we have that
		\begin{equation}
			\mathbf{J} = \lim _{n \rightarrow \infty} \frac{1}{n} \sum_{i=1}^{n} \mathbf{V}_i^T \mathsf{Cov} \left( \mathbb{E} \left [ \mathbf{T}_{X_i}|B_i\right ] \right ) \mathbf{V}_i,
		\end{equation}
		which needs to be positive definite with finite determinant. This is ensured by condition~3 in Theorem~\ref{thm:MainThm}.
		\item \label{con:2.3_multi}
  This condition is~\cite[condition II(iii)]{bradley1962} and, as the two above, it is needed to apply the weak law of large numbers for independent random variables~\cite[p.174]{cramer1947problems}.
		In order to check this condition  we first need to find the following set for all $i \in [1:n]$
		\begin{equation*}
			D_{5 i}\! = \!\left\{b \in \{-1,1\}: H_{irst}\left(b\right)>n  \right\}.
		\end{equation*}
		As discussed above in~\ref{con:1.2 multi},
		whenever condition~\ref{con:1 Thm} and condition~\ref{con:2 Thm} in Theorem~\ref{thm:MainThm} hold, we have that $H_{i r s t}$ is finite for all $i \in [1:n]$. Thus, for a sufficiently large $n$, there will be no $b \in\{-1,1\}$ such that $H_{i r s t}\left(b\right)>n$ for all $i \in [1:n]$. Hence, for large $n$, the set $D_{5 i}$ is empty, i.e.,
		\begin{equation}
			D_{5 i}=\varnothing \quad \text{ for all } \quad i \in[1:n].
		\end{equation}
		Hence, 
		\begin{equation}
			\sum_{i=1}^n \sum_{b \in D_{5 i}} P_{B_i}\left(b; \boldsymbol{\theta} \right)=0.
		\end{equation}
		This satisfies the condition $\sum_{i=1}^n \sum_{b \in D_{5 i}} P_{B_i}\left(b;\boldsymbol{\theta} \right)=o(1)$~\cite[condition II(iii)]{bradley1962}.
		Similarly, we can find the set $D_{6 i}$ for all $i \in [1:n]$ as follows,		
		\begin{equation}
			D_{6 i} =  \left\{b \in \{-1,1\}: H_{i r s t}\left(b\right) < n \right\}.
		\end{equation}
		Using the same reasoning as above we have that $D_{6i} = \{-1,1\}$ for a sufficiently large $n$.
		Now, we need to check the following condition~\cite[condition II(iii)]{bradley1962},
		\begin{equation}
			\sum_{i=1}^{n} \sum_{b \in D_{6 i}} H_{i r s t}^{2}\left(b\right) P_{B_i}(b;\boldsymbol{\theta}) = o\left(n^{2}\right).
		\end{equation}
		From~\ref{con:1.2 multi} above, we have that $H_{i r s t}\left(b\right) = K$ for all $(r, s, t) \in [1:k]^3$, $b \in \{-1,1\}$, and $i \in [1:n]$. Hence,
		\begin{align}
			& \sum_{i=1}^{n} \sum_{b \in D_{6 i}} H_{i r s t}^{2}\left(b\right) P_{B_i}(b;\boldsymbol{\theta}) \notag
			\\ & \leq \sum_{i=1}^{n} K^2 \sum_{b \in D_{6 i}}  P_{B_i}(b;\boldsymbol{\theta}) = K^2 n = o(n^2),
		\end{align}
		and
		\begin{equation}
			\frac{1}{n} \sum_{i=1}^{n} M_{i} = K ,
		\end{equation}	
		which is a finite positive constant.
  \begin{figure*}
  \begin{align}
  \label{eq:Con6GaussianAsy}
			& \left (\sum_{r=1}^{k}\left(\frac{\partial \log P_{B_i}(b; \boldsymbol{\theta})}{ \partial \theta_{r}}\right)^{2} \right)^{\frac{1}{2}} =\left(\sum_{r=1}^{k}\left(\sum_{j=1}^d \left(\mathbb{E}\left[\mathbf{T}_{X_i}(j)|B_i= b\right ]-\mathbb{E}\left[\mathbf{T}_{X_i}(j) \right ] \right ) \mathbf{V}_i(j,r)\right)^{2} \right)^{\frac{1}{2}}
		\end{align}
     \hrule
  \end{figure*}
		\item 
  This condition is~\cite[condition III]{bradley1962} and it is needed to ensure asymptotic normality. In particular, we need
		\begin{equation} \label{eq:con 3.1Multi}
			\lim _{n \rightarrow \infty}\!\frac{1}{n} \!\sum_{i=1}^{n}\sum_{b \in D_{7 i}} \sum_{r=1}^k\!\left(\!\frac{\partial \log P_{B_i}(b;\boldsymbol{\theta})}{\partial \theta_r}\!\right)^{2}\!\!\! P_{B_i}(b;\boldsymbol{\theta})\!=\!0,
		\end{equation}	
		where $D_{7 i}$ for all $i \in [1:n]$ is defined as follows,
		\begin{equation*}
			D_{7 i}\!=\!\!\Big\{\!b\!\in\!\{\!-1,1\!\}\!:\! \Big[\!\sum_{r=1}^{k}\!\Big(\frac{\partial \log P_{B_i}(b; \boldsymbol{\theta})}{ \partial \theta_{r}}\!\Big)^{2} \Big]^{\frac{1}{2}}\!\!>\!\epsilon \sqrt{n} \Big\}\!,
		\end{equation*}
		for every $\epsilon>0$.	By using Lemma~\ref{lem:Score}, we have that~\eqref{eq:Con6GaussianAsy} at the top of the next page holds.
		
		Moreover, as shown above in~\ref{con:1.1_multi}, whenever condition~\ref{con:1 Thm} and condition~\ref{con:2 Thm} in Theorem~\ref{thm:MainThm} hold,
		the derivative $ \partial \log P_{B_i} / \partial \theta_r$ is finite for all $i \in [1:n]$ and $r \in [1:k]$. 
		Thus, for all $i \in [1:n]$ and $r \in [1:k],$ there exists some positive constant $C$ such that
		\begin{equation}
			\sum_{j=1}^d \left(\mathbb{E} \left [ \mathbf{T}_{X_i}(j)|B_i \!=\! b\right ] \!-\! \mathbb{E}\left[\mathbf{T}_{X_i}(j) \right ] \right ) \mathbf{V}_i(j,r) \!\leq\! \sqrt{C},
		\end{equation}
		which leads to
		\begin{equation}
			\Big (\sum_{r=1}^{k}\Big(\frac{\partial \log P_{B_i}(b; \boldsymbol{\theta})}{ \partial \theta_{r}}\Big)^{2} \Big)^{\frac{1}{2}}\leq \Big( \sum_{r=1}^k C\Big )^\frac{1}{2} = \sqrt{Ck}.
		\end{equation}
		Thus, for a sufficiently large $n$, there will be no $b \in\{-1,1\}$ such that 
		\begin{equation}
			\Big(\sum_{r=1}^{k}\Big(\partial \log P_{B_i}(b) / \partial \theta_{r}\Big)^{2}\Big)^{\frac{1}{2}}>\epsilon \sqrt{n},
		\end{equation}
		for any $\epsilon>0$. Hence, for large $n$, the set $D_{7 i}$ will be empty for any $\epsilon>0$, i.e.,
		\begin{equation}
			D_{7i}=\varnothing, \quad \text{ for all } \quad  i \in[1:n].
		\end{equation}
	\end{enumerate}
    Thus, the condition in~\eqref{eq:con 3.1Multi} will be satisfied. This concludes the proof of Theorem~\ref{thm:MainThm}.
    \section{Derivation of the FIM for the Gaussian Distribution}
  	\label{app:FIMGaussian}
  	\subsection{FIM in Case~1} \label{app:FIM in Case.1}
    To compute the FIM in Proposition~\ref{prop:FIM}, 
    we need to derive the following variance,
	\begin{align} \label{eq:unknownmean}
		& \mathsf{Var} \left(\mathbb{E}[\mathrm{T}_{X_i} | B_i] \right) \notag 
   = \mathbb{E}\left[\left( \mathbb{E}[X_i | B_i]\right)^2\right] - \left(\mathbb{E}[X_i]\right)^2 \notag
		\\ & = \left(\mathbb{E}\left[X_i | B_i=1\right]\right)^2 \Pr(B_i=1)   \notag
            \\ & \quad +\left(\mathbb{E}\left[X_i | B_i=-1\right]\right)^2 \Pr(B_i=-1) - \left(\mathbb{E}[X_i]\right)^2 \notag
        \\ & =   \sigma^4 \frac{p_{X_i}^2(\tau_i)}{F_{X_i}(\tau_i) (1-F_{X_i}(\tau_i))},
	\end{align}
	where the last equality follows since
     \begin{equation}
    \mathbb{E}[X_i|B_i =b] = \mu_i -b \sigma^2 \frac{p_{X_i}(\tau_i)}{\Pr(B_i=b)}, \ b \in \{-1,1\}.
    \label{eq:CondExpectInterval}
     \end{equation}
By substituting~\eqref{eq:unknownmean} inside~\eqref{eq:FIM_censored} with $v_i = \frac{w_i}{\sigma^2}$, we obtain~\eqref{eq:FIMGaussCase1}.

  	\subsection{FIM in Case~2} \label{FIM in Case.2}
    To compute the FIM in Proposition~\ref{prop:FIM}, 
    we need to derive the following variance,
	\begin{align} \label{eq:unknowvariance}
		& \mathsf{Var} \left(\mathbb{E}[\mathrm{T}_{X_i} | B_i] \right) = \mathbb{E}\left[\left( \mathbb{E}[\mathrm{T}_{X_i} | B_i]\right)^2\right] - \left(\mathbb{E}[\mathrm{T}_{X_i}]\right)^2 \notag
  \\& =\mathbb{E}\left[\left( \mathbb{E}[(X_i-\mu_i)^2 | B_i]\right)^2\right] - \left(\mathbb{E}[(X_i-\mu_i)^2]\right)^2 \notag
		\\ & = \left(\mathbb{E}\left[(X_i-\mu_i)^2 | B_i=1\right]\right)^2 \Pr(B_i=1) \notag
        \\ & \quad + \left(\mathbb{E}\left[(X_i-\mu_i)^2 | B_i=-1\right]\right)^2 \Pr(B_i=-1) \notag
  \\& \quad- \left(\mathbb{E}[(X_i-\mu_i)^2]\right)^2 \notag
	\\& = \sigma^4 \left(\tau_i - \mu_i\right)^2\frac{p_{X_i}^2(\tau_i)}{F_{X_i}(\tau_i)(1-F_{X_i}(\tau_i))},
	\end{align}
	where the last equality follows from~\eqref{eq:CondExpectInterval} and the fact that for $\ b \in \{-1,1\}$ we have that
 \begin{equation}
 \label{eq:CondSecMomBound}
 \mathbb{E}[X_i^2|B_i =b] = \sigma^2 + \mu_i^2 -b \sigma^2  \frac{p_{X_i}(\tau_i)}{\Pr(B_i=b)} \left( \tau_i + \mu_i\right ).
 \end{equation}
Thus, we have that
  \begin{align}
  \label{eq:unknownvariance}
 & \mathbb{E}[(X_i - \mu_i)^2|B_i =b] \notag
 \\&= \sigma^2 + \mu_i^2 -b \sigma^2  \frac{p_{X_i}(\tau_i)}{\Pr(B_i=b)} \left( \tau_i + \mu_i\right )   \notag
 \\ & \quad - 2 \mu_i\left(\mu_i - b \sigma^2  \frac{p_{X_i}(\tau_i)}{\Pr(B_i=b)}\right)+ \mu_i^2 \notag
 \\ & = \sigma^2 -b \sigma^2  \frac{p_{X_i}(\tau_i)}{\Pr(B_i=b)} \left( \tau_i - \mu_i\right ).
 \end{align}
By substituting~\eqref{eq:unknownvariance} inside~\eqref{eq:FIM_censored} with $v_i = -1/2$, we obtain~\eqref{eq:FIMGaussCase2}. 

   \subsection{FIM in Case~3} \label{FIM in Case.3}
   To compute the FIM in Proposition~\ref{prop:FIM}, 
    we need to derive $\mathsf{Cov} \left(\mathbb{E}[\mathbf{T}_{X_i} | B_i] \right)$.
    In what follows, we let $\mathbf{T}_{X_i}(j), j \in [1:2]$ denote the $j$-th component of $\mathbf{T}_{X_i}$.
    We start by noting that from Case~1 in Appendix~\ref{app:FIM in Case.1}, we have that
	\begin{equation} \label{eq:var1}
		\mathsf{Var} \left(\mathbb{E}[\mathbf{T}_{X_i}(1) | B_i] \right) = \sigma^4 \frac{p_{X_i}^2(\tau_i)}{F_{X_i}(\tau_i) (1-F_{X_i}(\tau_i))}.
	\end{equation} 
 Similarly, from Case~2 (see~\eqref{eq:unknowvariance}) in Appendix~\ref{FIM in Case.2}, we obtain
	\begin{equation} \label{eq:var2}
		\mathsf{Var} \left(\mathbb{E}[\mathbf{T}_{X_i}(2) | B_i] \right)  =\sigma^4 \left(\tau_i + \mu_i\right)^2\frac{p_{X_i}^2(\tau_i)}{F_{X_i}(\tau_i)(1-F_{X_i}(\tau_i))},
	\end{equation}
 where $\mu_i = w_i\alpha$.
	Finally, we have that
	\begin{align} \label{eq:cov}
		& \mathsf{Cov}\left(\mathbb{E}[\mathbf{T}_{X_i}(1) | B_i],\mathbb{E}[\mathbf{T}_{X_i}(2) | B_i] \right) \notag
  \\& = \mathsf{Cov}\left(\mathbb{E}[X_i | B_i],\mathbb{E}[X_i^2 | B_i] \right) \notag
  \\& = \mathbb{E} \left [\mathbb{E}[X_i | B_i]\mathbb{E}[X_i^2 | B_i] \right ] - \mathbb{E} \left [ \mathbb{E}[X_i | B_i]\right ]\mathbb{E} \left [ \mathbb{E}[X_i^2 | B_i]\right ] \notag
		\\ & = \mathbb{E}[X_i| B_i=1] \mathbb{E}[X_i^2 | B_i=1]\Pr\left(B_i=1\right) \notag
    \\ & \quad +  \mathbb{E}[X_i | B_i=-1] \mathbb{E}[X_i^2 | B_i=-1]\Pr\left(B_i = -1\right) \notag	  
    \\& \quad-  \mathbb{E}[X_i] \mathbb{E}[X_i^2] \notag
  \\& = \mu_i \sigma^2 + \mu_i^3 + \sigma^4 (\tau_i+\mu_i)  \frac{p_{X_i}^2(\tau_i)}{F_{X_i}(\tau_i)(1-F_{X_i}(\tau_i))}  \notag
  \\ & \quad  - \mu_i(\sigma^2 + \mu_i^2) \notag
  \\& = \sigma^4 (\tau_i+\mu_i) \frac{p_{X_i}^2(\tau_i)}{F_{X_i}(\tau_i)(1-F_{X_i}(\tau_i))},
	\end{align}
 where we have used~\eqref{eq:CondExpectInterval} and~\eqref{eq:CondSecMomBound} to find $\mathbb{E}[X_i|B_i =b]\mathbb{E}[X_i^2|B_i =b]$ for $b \in \{-1,1\}$.
By substituting~\eqref{eq:cov} inside~\eqref{eq:FIM_censored} with $\mathbf{V}_i = \begin{bmatrix}
			w_i &0
			\\ 0 & -\frac{1}{2}
		\end{bmatrix}$, we obtain~\eqref{eq:FIM case3}.

    \section{Proof of Proposition~2}
	\label{app:Proofprop2}
%
%
%
%
We need to make sure that 
\begin{equation}
\begin{bmatrix} u & v \end{bmatrix} \mathbf{J} \begin{bmatrix} u & v \end{bmatrix}^T >0
\end{equation}
for all vectors $\begin{bmatrix} u & v \end{bmatrix}^T \neq \mathbf{0}_2$. By using $\mathbf{J}_n$ in~\eqref{eq:FIM case3}, the above condition can be equivalently written as
\begin{equation}
\label{eq:Cond3PosCase}
\lim_{n \to \infty} \frac{1}{n} \sum_{i=1}^n c_i \left(w_i u - v \frac{(\tau_i+w_i \alpha)}{2} \right )^2 >0.
\end{equation} 
It is now a simple exercise to show that if the  $\tau_i$'s, with $i \in [1:n]$, are chosen i.i.d. from some absolutely continuous distribution, then~\eqref{eq:Cond3PosCase} is satisfied almost surely.
This concludes the proof of Proposition~\ref{prop:Case3GaussCond}.

    \section{Derivation of the FIM for the Poisson Distribution} \label{app:FIMpoisson}
The probability mass function of $X_i$ in~\eqref{eq:systemModel2} is given by
  \begin{align}
        p_{X_i}(x)
    = \frac{\exp (v_i \theta x - \exp(v_i\theta))}{x!}. 
    \end{align}
With reference to~\eqref{eq:PDFExpFamily}, we have that $\mathrm{T}_x =  x$ and
\begin{subequations}
\begin{align*}
&h(x) =\frac{1}{x\ !},
\\& \eta_{i,\theta} = v_i\theta, 
\\& \phi(\eta_{i,\theta}) = \exp(v_i\theta) = \exp(\eta_{i,\theta}).
\end{align*}
\end{subequations}
Now, to compute the FIM in Proposition~\ref{prop:FIM} for the model in~\eqref{eq:systemModel2}, 
    we need to derive the following variance,   
	\begin{align} \label{eq:unknowntheta}
		& \mathsf{Var} \left(\mathbb{E}[\mathrm{T}_{X_i} | B_i] \right)  
   = \mathbb{E}\left[\left( \mathbb{E}[X_i | B_i]\right)^2\right] - \left(\mathbb{E}[X_i]\right)^2 \notag
\\ & = \left(\mathbb{E}\left[X_i | B_i=1\right]\right)^2 \Pr(B_i=1)   \notag
            \\ & \quad +\left(\mathbb{E}\left[X_i | B_i=-1\right]\right)^2 \Pr(B_i=-1) - \left(\mathbb{E}[X_i]\right)^2 \notag
  \\ &  = \exp (2v_i\theta) \frac{\left(F_{X_i}( \tau_i ) - F_{X_i}( \tau_i -1)\right)^2}{F_{X_i}( \tau_i )\left(1-F_{X_i}( \tau_i )\right)},
    \end{align}
where the last equality follows since for $b \in \{-1,1\}$ it holds that
     \begin{equation}
    \mathbb{E}[X_i|B_i\!=\!b]\!=\! \frac{{\rm{e}}^{\!v_i\theta}\!\left(1\!\!-\!F_{X_i}( \tau_i \!\!-\!\!1)\right)^{\frac{1\!-\!b}{2}}\!\left(F_{X_i}( \tau_i \!-\!1)\right)^{\frac{1\!+\!b}{2}}}{\left(1\!-\!F_{X_i}( \tau_i )\right)^{\frac{1-b}{2}}\left(F_{X_i}( \tau_i )\right)^{\frac{1+b}{2}}}.
    \label{eq:CondExpectInterval_poisson}
     \end{equation}
By substituting~\eqref{eq:unknowntheta} inside~\eqref{eq:FIM_censored}, we obtain~\eqref{eq:FIMPoisson}. 
\bibliographystyle{IEEEtran}
\bibliography{refs}
\end{document}